\documentclass[11pt]{amsart}
\usepackage{amsfonts,latexsym}
\usepackage{amsmath}
\usepackage{amscd}
\usepackage{float,amsmath,amssymb,mathrsfs,bm,multirow,graphics}
\newtheorem{theorem}{Theorem}
\newtheorem{proposition}[theorem]{Proposition}

\newtheorem{remark}[theorem]{Remark}

\newcommand{\proofbegin}{\noindent{\it Proof.\,\,}}
\newcommand{\proofend}{\hfill$\Box$\bigskip}

\newcommand{\pot}{{\bf Proof of Theorem \ref{th:lwp}}}
\newcommand{\pott}{{\bf Proof of theorem \ref{th:analytic}}}

\title{The periodic Cauchy problem for Novikov's equation}
\author{Feride T\i\u{g}lay 
}
\address{Fields Institute, 222 College Street, 2nd Floor, Toronto, Ontario M5T 3J1, Canada}
\email{ftiglay@fields.utoronto.ca} 
\begin{document}
\begin{abstract} We study the periodic Cauchy problem for an integrable equation with cubic nonlinearities introduced by V. Novikov. We show the local well-posedness of the problem in Sobolev spaces and existence and uniqueness of solutions for all time using orbit invariants. Furthermore we prove a Cauchy-Kowalevski type theorem for this equation, that establishes the existence and uniqueness of real analytic solutions. 
\end{abstract}  
%
%
\maketitle
\section{Introduction}
Recently the integrable equation with cubic nonlinearities
\begin{equation}
\label{eq:nov1}
u_t-u_{xxt}+4u^2u_x-3uu_xu_{xx}-u^2u_{xxx}=0
\end{equation}
derived by V. Novikov in \cite{N} has attracted some attention in the litterature \cite{HLS, HW}. We study the periodic Cauchy problem for this equation for Sobolev class and real analytic data. We prove local (in time) well-posedness in Sobolev spaces and existence and uniqueness of Sobolev class solutions for all time provided that the initial data of the same class satisfies a sign condition. Furthermore we prove a Cauchy-Kowalevski type theorem for this equation using a contraction argument on a decreasing scale of Banach spaces.

Note that one can write Novikov's equation \eqref{eq:nov1} in the form 
\begin{equation}
\label{eq:nov2}
m_t+m_xu^2+3muu_x=0 \quad \mbox{where} \quad m=u-u_{xx}.
\end{equation}

Like the Camassa-Holm and Degasperis-Procesi equations, Novikov's equation \eqref{eq:nov1} has Lax pair representations and admits peakon solutions, but it has nonlinear terms that are cubic, rather than quadratic. A Lax representation for this equation is introduced by Novikov in \cite{N} in the form
\[ \psi_{xxx}=\psi_x + \lambda m^2 \psi+2\frac{m_x}{m}\psi_{xx}+ \frac{mm_{xx}-2m_x^2}{m^2}\psi_x,
\]
\[ \psi_t= \frac{u}{\lambda m} \psi_{xx}-\frac{mu_x+um_x}{m^2}\psi_x-u^2\psi_x.
\]

Another representation of this equation as a matrix Lax pair 
is given in \cite{HW} and is shown to be related to a negative flow in the Sawada-Kotera hierarchy. In the same article a bi-Hamiltonian structure is provided without proof in the form 
$m_t=B_1(\delta H_1/\delta m)=B_2(\delta H_2/\delta m)$ 
where the operators are $B_1=-2(3m\partial_x+2m_x)(4\partial_x-\partial_x^3)^{-1}(3m\partial_x+m_x)$ and $B_2=(1-\partial_x^2)\frac{1}{m}\partial_x \frac{1}{m}(1-\partial_x^2)$
and the Hamiltonians are $H_1=\frac{1}{3}\int (m^{-8/3}m_x^2+9m^{-2/3})dx $ and $ H_2= \frac{1}{8}\int (u^4+2u^2u_x^2-\frac{1}{3}u_x^4)dx$.

Another interesting property of Novikov's equation \eqref{eq:nov1} that is common with Camassa-Holm and Degasperis-Procesi equations is that it has nonsmooth soliton solutions with multiple peaks (multipeakons). Multipeakons for \eqref{eq:nov1} are explicitly computed in \cite{HLS} using scattering theory.   

In section 2 we prove local well-posedness (existence, uniqueness and continuous dependence on initial data for a short time), persistence of solutions and global existence and uniqueness of solutions for the periodic Cauchy problem for Novikov's equation. In order to prove local well-posedness  we use a method that can be traced back to the work of D.G. Ebin and J. Marsden \cite{EM}. Following the elegant geometric framework given by V. Arnold in \cite{A}, they develop in \cite{EM} the analytic tools and prove sharp local well-posedness results in Sobolev and H\"{o}lder spaces for the Dirichlet problem for Euler equations of ideal hydrodynamics. This method has also been implemented to prove local well-posedness of the periodic Cauchy problems for nonlinear evolution equations such as Camassa-Holm and Hunter-Saxton equations (see \cite{HM2}, \cite{Mis02} and \cite {Ti1}). Here we use this method for Novikov's equation \eqref{eq:nov1} and prove a local well-posedness theorem in Sobolev spaces.      

Our global existence and uniqueness theorem \ref{th:gwp} uses orbit invariants to establish that the solutions persist for all time if a sign condition holds. We refer to \cite{CM} for a similar result for Camassa-Holm equation on the real line. A detailed discussion of these orbit invariants is given in \cite{TV}. 

In section 3 we prove the analytic regularity (i.e., existence and uniqueness of analytic solutions for analytic initial data) of the Cauchy problem for Novikov's equation. It is well known that the solutions to the Hunter-Saxton and Camassa-Holm equations are analytic in both space and time variables for a short time (see \cite{Ti1} and \cite{HM3} respectively).
In contrast the solutions of the Korteweg-De Vries equation are analytic in the space variable for all time \cite{Tru} but are not analytic in the time variable \cite{KaM}. Theorem \ref{th:analytic} establishes that solutions of the Novikov's equation, like Hunter-Saxton and Camassa-Holm equations, are analytic in both space and time variables.

Our approach in proving theorem \ref{th:analytic} is to use a contraction argument on an appropriate scale of Banach spaces. The general framework for this existence theorem has been developed as an abstract Cauchy-Kowalevski theorem by L.V. Ovsjannikov \cite{Ovs1, Ovs2}, F. Treves \cite{Tre}, L. Nirenberg \cite{Nir}, T. Nishida \cite{Nis} and M.S. Baouendi and C. Goulaouic \cite{BG} among others and subsequently applied to the Euler and Navier-Stokes equations.

\section{Local well-posedness and existence of global solutions}

We formulate the periodic Cauchy problem for \eqref{eq:nov1} as follows using the inverse $\Lambda^{-2}$ of the elliptic operator $\Lambda^2=1-\partial_x^2$:
\begin{eqnarray}
& & u_t+u^2u_x=-\Lambda^{-2}\left( 3u^2u_x+2u_x^3+3uu_xu_{xx}\right), \quad \quad x\in{\Bbb T}={\Bbb R}/{\Bbb Z}, \ t\in{\Bbb R}, \label{eq:nloc} \\
& & u(0,x)=u_0(x). \label{eq:data}
\end{eqnarray}

In our estimates we use some properties of Sobolev class functions. We summarize them here for the convenience of the reader. Note that all bounds are up to a constant that may depend on $H^s$ norms of the Sobolev class diffeomorphisms $\eta$ and $\eta^{-1}$. 

Let $s>3/2$. We denote the space of circle diffeomorphisms of Sobolev class $H^s$ by $\mathrm{Diff}^s(S^1)$. If $u\in H^s$ then the composition map $\eta\mapsto u\circ\eta$ from $\mbox{Diff}^s$ to $H^s$ and the inversion map $\eta\mapsto\eta^{-1}$ on $\mbox{Diff}^s$ are continuous. Furthermore, for any $\eta\in\mbox{Diff}^s(S^1)$, \footnote{From this point on, $\lesssim$ denotes an inequality up to a constant that depends only on $\|\eta\|_{H^s}$ and $\|\eta^{-1}\|_{H^s}$.}
\begin{equation}\label{eq:Hscomp}
\|u\circ\eta\|_{H^s} \lesssim 
(1+\|\eta\|_{H^s}^s)\|u\|_{H^s}.
\end{equation}

Another tool we use is a commutator estimate from \cite{KP}. For $s>0$ and $\Lambda^s=(1-\partial_x^2)^{s/2}$, if $u,v\in H^s({\Bbb T})$ then
\begin{equation}\label{eq:commutator}
\| [\Lambda^s,u]v\|_{L^2}\lesssim \|\partial_x u\|_{\infty} \|\Lambda^{s-1}v\|_{L^2}+ \|\Lambda^s u\|_{L^2} \|v\|_{\infty}.
\end{equation}

\begin{theorem}[Local well-posedness] \label{th:lwp}
For $s>5/2$ and $u_0\in H^s({\Bbb T})$ there is a $T>0$ and a unique solution
\[ u\in C((0,T), H^s({\Bbb T}))\cap C^1((0,T), H^{s-1}({\Bbb T}))
\]
to the problem \eqref{eq:nloc}-\eqref{eq:data} that depends continuously on initial data.
\end{theorem}
Our strategy will be to reformulate \eqref{eq:nloc}-\eqref{eq:data} 
as an initial value problem on the space of circle diffeomorphisms $\mathrm{Diff}^s(S^1)$ of Sobolev class $H^s$. 
It is well known that whenever $s>3/2$ this space is a smooth Hilbert manifold and a topological group. 
We will then show that the reformulated problem can be solved 
on $\mathrm{Diff}^s(S^1)$ by standard ODE techniques. 

It is convenient to introduce the notation 
$A_{\xi}=R_{\xi}\circ A \circ R_{\xi^{-1}}$
for the conjugation of an operator $A$ on $H^s(S^1)$ by 
a diffeomorphism $\xi \in \mathrm{Diff}^s(S^1)$, for instance $\partial_{x_{\xi}}f$ means $(\partial_x (f\circ\xi^{-1}))\circ\xi$. 

Let $u=u(t,x)$ be a solution of \eqref{eq:nloc} with initial data $u_0$. 
Then the flow $t \to \xi(t, x)$ associated to $u^2$, i.e. the solution of 
the initial value 
problem\footnote{Here ``dot'' indicates differentiation in $t$ variable.} 
$$\dot\eta(t,x)=u^2(t,\eta(t,x)), \quad \eta(0,x)=x$$
is (at least for a short time) a smooth curve in the space of diffeomorphisms 
starting from the identity $\mbox{id} \in \mathrm{Diff}^s(S^1)$. On the other hand a solution $(\eta, \zeta)$ of the initial value problem
\begin{eqnarray}
& & \dot{\eta}=\zeta^2, \\
& & \dot{\zeta}=-\Lambda^{-2}_{\eta}\left\{ 3\zeta^2 \partial_{x_{\eta}} \zeta+2(\partial_{x_{\eta}} \zeta)^3+3\zeta\partial_{x_{\eta}} \zeta \partial_{x_{\eta}}^2 \zeta \right\}= F(\eta,\zeta), \\
& & \eta(0,x)=x, \quad \zeta(0,x)=u_0(x) 
\end{eqnarray}
determines a solution $u(t,x)=\zeta(t,\eta^{-1}(t,x))$ of the problem \eqref{eq:nloc}-\eqref{eq:data}.

In the proof of theorem \ref{th:lwp} we make repeated use of the estimates in \eqref{eq:Hscomp}, \eqref{eq:commutator} and Sobolev embedding theorems.
\medskip

\pot
Note that it suffices to prove that the map $(\eta, \zeta)\rightarrow F(\eta, \zeta)$ is Fr\'{e}chet differentiable. We show that $F(\eta,\zeta)$  maps $\mathrm{Diff}^s(S^1)\times H^s({\Bbb T})$ into $H^s({\Bbb T})$ and that its directional derivatives $\partial_{\eta}F_{(\eta,\zeta)}$ and $\partial_{\zeta}F_{(\eta,\zeta)}$ define bounded linear maps that are continuous in both $\eta$ and $\zeta$. Then Fr\'{e}chet differentiability of $(\eta, \zeta)\rightarrow F(\eta, \zeta)$ follows.

Our first estimate establishes the boundedness of the map $(\eta, \zeta)\rightarrow F(\eta, \zeta)$: We use the ring property of Sobolev spaces with the estimate \eqref{eq:Hscomp} to obtain
\begin{eqnarray}
&  \|F\|_{H^s} & \lesssim \| (\zeta\circ\eta^{-1})^2 \partial_x(\zeta\circ\eta^{-1})\|_{H^{s-2}}+\| (\partial_x(\zeta\circ\eta^{-1}))^3 \|_{H^{s-2}} \\
&  & \quad + \| (\zeta\circ\eta^{-1})\partial_x(\zeta\circ\eta^{-1})\partial_x^2 (\zeta\circ\eta^{-1})\|_{H^{s-2}}\\
& & \lesssim \| \zeta\circ\eta^{-1}\|_{\infty}^2 \| \zeta\circ\eta^{-1}\|_{H^{s-1}}+ \| \zeta\circ\eta^{-1}\|_{H^{s-1}}^2 \| \zeta\circ\eta^{-1}\|_{C^1} \\
& & \quad + \| \zeta\circ\eta^{-1}\|_{\infty} \| \zeta\circ\eta^{-1}\|_{C^1} \| \partial_x^2(\zeta\circ\eta^{-1})\|_{H^{s-2}}.
\end{eqnarray}
Then the same tools imply
\[ \|F\|_{H^s} \leq C_{\eta} \| \zeta \|_{H^s}
\]
where $C_{\eta}$ depends only on the $H^s$ norms of $\eta$ and $\eta^{-1}$.
 
{\em Directional derivatives.} The derivative of $F$ in the direction of $\zeta$ is 
\begin{equation}\label{eq:der_zeta}
\partial_{\zeta}F_{(\eta,\zeta)}(X)= -\Lambda_{\eta}^{-2}\big\{ 6\partial_{x_{\eta}} \zeta \partial_{x_{\eta}}(X\zeta)+ 3\zeta\partial_{x_{\eta}}(\zeta\partial_{x_{\eta}} X)+ 3\partial_{x_{\eta}}(X \zeta)\partial_{x_{\eta}}^2\zeta\big \}.
\end{equation}
The same type of argument we gave above for the boundedness of the map defined by $F$ applies to $\partial_{\zeta}F_{(\eta,\zeta)}(X)$ as well, hence $X\rightarrow \partial_{\zeta}F_{(\eta,\zeta)}(X)$ is a bounded linear map on $H^s({\Bbb T})$. 
We compute $\partial_{\eta}F_{(\eta,\zeta)}(X)$ in steps to simplify the notation. First observe that we can write $\partial_{\eta}F_{(\eta,\zeta)}(X)$ as a sum
\begin{equation}\label{eq:der_eta}
\partial_{\eta}F_{(\eta,\zeta)}(X)=\frac{d}{ds}\Big|_{s=0}\big( \Lambda^{-2}_{\eta_s}\big)(W)
+ \Lambda^{-2}_{\eta}\Big(\frac{d}{ds}\Big|_{s=0}W(\eta_s,\zeta)\Big)\circ\eta
\end{equation}
where 
\[ W(\eta,\zeta)= 3\zeta^2 \partial_{x_{\eta}} \zeta +2 (\partial_{x_{\eta}} \zeta)^3 + 3\zeta \partial_{x_{\eta}} \zeta \partial_{x_{\eta}}^2 \zeta.
\]
The following formulas for the directional derivatives of the operators $\Lambda^{-2}_{\eta}$ and $\partial_{x_{\eta}}$ are straightforward to compute: Let $\eta_s|_{s=0}=\mbox{id}$ and $(d\eta_s/ds)|_{s=0}=X$. Then
\begin{equation}\label{eq:ddlambda}
\frac{d}{ds}\Big|_{s=0}\Lambda^{-2}_{\eta_s}(V)=\big[X, \Lambda^{-2}_{\eta}\big]\partial_{x_{\eta}} V,
\end{equation}
where $[\ . \ , \ . \ ]$ denotes the commutator, and
\begin{equation}\label{eq:ddderivative}
\frac{d}{ds}\Big|_{s=0}\partial_{x_{\eta}}(V)=-\partial_{x_{\eta}} V \partial_{x_{\eta}} X.
\end{equation}
Using \eqref{eq:ddderivative} we find the directional derivative of $W$ in $\eta$ to be
\begin{equation}\label{eq:ddW}
\frac{d}{ds}\Big|_{s=0}W(\eta_s,\zeta)= -3 \zeta^2 \partial_{x_{\eta}} \zeta \partial_{x_{\eta}} X -3\partial_{x_{\eta}} \big(\zeta (\partial_{x_{\eta}}\zeta)^2 \partial_{x_{\eta}} X\big)-3 \partial_{x_{\eta}} X \partial_{x_{\eta}} \zeta \partial_{x_{\eta}} \big(\zeta\partial_{x_{\eta}}\zeta\big). 
\end{equation}
Note that the cubic nonlinearities are not creating any extra difficulty here. Combining \eqref{eq:der_eta}, \eqref{eq:ddlambda} and \eqref{eq:ddW} we obtain
\begin{equation}
\partial_{\eta}F_{(\eta,\zeta)}(X)= [X, \Lambda^{-2}_{\eta}\partial_{x_{\eta}}] W -3 \Lambda^{-2}_{\eta}\partial_{x_{\eta}} \big( \zeta (\partial_{x_{\eta}} \zeta)^2 \partial_{x_{\eta}} X \big)-3 \Lambda^{-2}_{\eta}\big( (\partial_{x_{\eta}}\zeta)^3 \partial_{x_{\eta}} X \big).
\end{equation}
Observe that, by commutator and product estimates for Sobolev spaces,  $X \mapsto \partial_{\eta}F_{(\eta,\zeta)}(X)$ is a bounded linear map on $H^s$. 

Next we prove the continuity of the map $(\eta, \zeta) \mapsto \partial_{\eta}F_{(\eta,\zeta)}$ from $\mathrm{Diff}^s\times H^s$ to $L(H^s, H^s)$.
Note that it suffices to estimate 
\begin{equation}\label{eq:cont_deta}
\| \partial_{\eta}F_{(\eta,\zeta)}(X) - \partial_{\eta}F_{(\mbox{id},\zeta)}(X)\|_{H^s}
+ \| \partial_{\eta}F_{(\eta,\zeta_1)}(X) - \partial_{\eta}F_{(\eta,\zeta_2)}(X)\|_{H^s}.
\end{equation}
The first summand in \eqref{eq:cont_deta} is clearly bounded by
\begin{eqnarray}
&  \|[X, \Lambda^{-2}_{\eta}\partial_{x_{\eta}}]W(\eta,\zeta)- [X,\Lambda^{-2}\partial_x] W(\mbox{id},\zeta)\|_{H^s} \label{eq:oneA}\\
& +\| \Lambda^{-2}_{\eta}\partial_{x_{\eta}} \big( \zeta (\partial_{x_{\eta}} \zeta)^2 \partial_{x_{\eta}} X \big)- \Lambda^{-2}\partial_x \big(\zeta (\partial_x \zeta)^2 \partial_x X\big)\|_{H^s} \label{eq:oneB}\\
&  + \| \Lambda^{-2}_{\eta}\big( (\partial_{x_{\eta}}\zeta)^3 \partial_{x_{\eta}} X \big)- \Lambda^{-2}\big( (\partial_x\zeta)^3 \partial_x X\|_{H^s}. \label{eq:oneC}
\end{eqnarray}
Adding and subtracting appropriate terms in the norm \eqref{eq:oneA} and using composition properties of Sobolev spaces (see \cite{BB} and \cite{Mis02} for instance) we obtain the following bound for \eqref{eq:oneA}:
\begin{eqnarray}
& \| \eta - \mbox{id}\|_{H^s} \| [X\circ\eta^{-1}, \Lambda^{-2}\partial_x](W(\eta,\zeta)\circ\eta^{-1})\|_{H^s} \label{eq:oneA1}\\
& + \| [X\circ\eta^{-1}, \Lambda^{-2}\partial_x](W(\eta,\zeta)\circ\eta^{-1}-W(\mbox{id},\zeta))\|_{H^s} \label{eq:oneA2}\\
& +\|[X\circ\eta - X, \Lambda^{-2}\partial_x](W(\mbox{id},\zeta))\|_{H^s} \label{eq:oneA3}
\end{eqnarray}
Note that \eqref{eq:oneA1} is easily bounded by $\| \eta-\mbox{id}\|_{H^s} \|X\|_{H^s} \| W(\eta, \zeta)\|_{H^{s-2}}$ using commutator estimates \cite{Taylor} and the above mentioned 
properties of Sobolev spaces. For \eqref{eq:oneA2} and \eqref{eq:oneA3} we proceed similarly and obtain the respective bounds $\|X\|_{H^s} \| W(\eta,\zeta)\circ\eta^{-1}-W(\mbox{id},\zeta)\|_{H^{s-2}}$ and $\|X\circ\eta - X\|_{H^s} \|W(\mbox{id},\zeta)\|_{H^{s-2}}$. Note that we can write $W$ as
\begin{equation}
W(\eta,\zeta)= \partial_{x_{\eta}} \big( \zeta^3+\frac{3}{2}\zeta(\partial_{x_{\eta}}\zeta)^2\big)+\frac{1}{2}(\partial_{x_{\eta}}\zeta)^3
\end{equation}
and decompose the commutators in \eqref{eq:oneA2} and \eqref{eq:oneA3} similarly. Note also that
$\| W \|_{H^{s-2}}$ is bounded by $\| \zeta\|_{H^s}^3$. Therefore, both \eqref{eq:oneA2} and \eqref{eq:oneA3} are estimated by $\|\eta -\mbox{id}\|_{H^s}\|X\|_{H^s} \|\zeta\|_{H^s}^3$. Therefore \eqref{eq:oneA} is estimated by  $\| \eta-\mbox{id}\|_{H^s} \|X\|_{H^s} \| \zeta\|_{H^s}^3$.

Adding and subtracting the appropriate terms to \eqref{eq:oneB} and using composition properties of Sobolev spaces lead to
\begin{eqnarray}
& \| \Lambda^{-2}_{\eta}\partial_{x_{\eta}} \big( \zeta (\partial_{x_{\eta}} \zeta)^2 \partial_{x_{\eta}} X \big)- \Lambda^{-2}\partial_x \big(\zeta (\partial_x \zeta)^2 \partial_x X\big)\|_{H^s} \nonumber \\
& \lesssim \|\eta -\mbox{id}\|_{H^s} \| \Lambda^{-2}\partial_x \Big( \zeta\circ\eta^{-1} \big(\partial_x(\zeta\circ\eta^{-1})\big)^2 \partial_x(X\circ\eta^{-1})\Big)\|_{H^s} \\
& + \| \zeta\circ\eta^{-1} \big(\partial_x(\zeta\circ\eta^{-1})\big)^2 \partial_x(X\circ\eta^{-1}) - \zeta (\partial_x \zeta)^2 X\|_{H^{s-1}}.
\end{eqnarray}
Then we bound \eqref{eq:oneB} by $\| \eta-\mbox{id}\|_{H^s} \|X\|_{H^s} \| \zeta\|_{H^s}^3$ using \eqref{eq:Hscomp} and composition and algebra properties of Sobolev spaces. 

We proceed similarly for \eqref{eq:oneC} to obtain
\begin{eqnarray}
& \| \Lambda^{-2}_{\eta}\big( (\partial_{x_{\eta}}\zeta)^3 \partial_{x_{\eta}} X \big)- \Lambda^{-2}\big( (\partial_x\zeta)^3 \partial_x X\|_{H^s} \nonumber \\
& \lesssim \| \eta-\mbox{id}\|_{H^s} \| \big(\partial_x (\zeta\circ\eta^{-1})\big)^3 \partial_x(X\circ\eta^{-1}) \|_{H^{s-2}} \\
& + \| \big(\partial_x (\zeta\circ\eta^{-1})\big)^3 \partial_x(X\circ\eta^{-1})-(\partial_x \zeta)^3 \partial_x X\|_{H^{s-2}}
\end{eqnarray}
from which follows the estimate $\| \eta-\mbox{id}\|_{H^s} \|X\|_{H^s} \| \zeta\|_{H^s}^3$ for \eqref{eq:oneC}.

Let us now consider the second summand in \eqref{eq:cont_deta} and combine the matching terms to estimate $\| \partial_{\eta}F_{(\eta,\zeta_1)}(X) - \partial_{\eta}F_{(\eta,\zeta_2)}(X)\|_{H^s}$ by
\begin{eqnarray}
& \lesssim \| [X, \Lambda^{-2}_{\eta}\partial_{x_{\eta}}] \big(W(\eta,\zeta_1)-W(\eta, \zeta_2)\big) \|_{H^s} \\
 & + \| \Lambda^{-2}_{\eta} \partial_{x_{\eta}} \big( \zeta_1 (\partial_{x_{\eta}}\zeta_1)^2 \partial_{x_{\eta}} X - \zeta_2 (\partial_{x_{\eta}}\zeta_2)^2 \partial_{x_{\eta}} X\big)\|_{H^s} \\
 & + \| \Lambda^{-2}_{\eta} \big( (\partial_{x_{\eta}} \zeta_1)^3 \partial_{x_{\eta}} X - (\partial_{x_{\eta}} \zeta_2)^3 \partial_{x_{\eta}} X \big)\|_{H^s}.
\end{eqnarray}
Using algebra property and Sobolev imbedding theorem along with the commutator estimate \eqref{eq:commutator} we have 
\begin{eqnarray}
& \| \partial_{\eta}F_{(\eta,\zeta_1)}(X) - \partial_{\eta}F_{(\eta,\zeta_2)}(X)\|_{H^s} \nonumber \\
& \lesssim \|X\|_{H^s} \| W(\eta,\zeta_1)-W(\eta, \zeta_2)\|_{H^{s-2}} \\
& + \|X\|_{H^s} \| \zeta_1 \circ\eta^{-1} \big( \partial_x(\zeta_1\circ\eta^{-1})\big)^2 - \zeta_2 \circ\eta^{-1} \big( \partial_x(\zeta_2\circ\eta^{-1})\big)^2\|_{H^{s-1}} \\
& + \|X\|_{H^s} \| \big( \partial_x(\zeta_1\circ\eta^{-1})\big)^3 - \big( \partial_x(\zeta_2\circ\eta^{-1})\big)^3 \|_{H^{s-2}}.
\end{eqnarray}
Adding and subtracting appropriate terms and using the same tools as above we bound $\| \partial_{\eta}F_{(\eta,\zeta_1)}(X) - \partial_{\eta}F_{(\eta,\zeta_2)}(X)\|_{H^s}$ by $\|X\|_{H^s} \|\zeta_1 -\zeta_2\|_{H^s}^3$. Therefore the map $(\eta, \zeta) \mapsto \partial_{\eta}F_{(\eta,\zeta)}$ from $\mathrm{Diff}^s\times H^s$ to $L(H^s,H^s)$ is continuous.

The continuity of $(\eta, \zeta) \mapsto \partial_{\zeta}F_{(\eta,\zeta)}$ from $\mathrm{Diff}^s\times H^s$ to $L(H^s, H^s)$ follows using the same tools except the commutator estimates. This completes the proof of theorem \ref{th:lwp}.
\proofend

Our next step towards proving existence of solutions for all time is a result stated in proposition \ref{prop:persistence} that establishes a condition which, if satisfied, guarantees that short time solutions persist for all time. This approach is in the same spirit as a  persistence result by J. T. Beale, T. Kato and A. Majda \cite{BKM} for solutions of the Euler equations for ideal hydrodynamics.    

\begin{proposition}[Persistence of solutions]\label{prop:persistence}
Let $s>5/2$ and $u\in C\big( [0,T], H^s({\Bbb T})\big)$ be a solution of the Cauchy problem \eqref{eq:nloc}-\eqref{eq:data}. If there exists a constant $K>0$ such that 
\[ \|u\|_{C^1}\leq K
\]
for all $t$ then the solution $u$ can be extended to a solution that persists for all time.
\end{proposition}
\proofbegin
The statement follows from Gronwall's inequality. We show that the differential inequality
\begin{equation}\label{eq:persgoal}
\frac{d}{dt}\|u\|_{H^s}^2 \lesssim \|u\|_{C^1}^2 \|u\|_{H^s}^2
\end{equation}
holds, hence $\|u\|_{H^s}$ stays bounded as long as $\|u\|_{C^1}$ is.

For a solution $u$ of \eqref{eq:nloc} we have
\begin{eqnarray}
&\frac{d}{dt}\|J_{\varepsilon}u\|_{H^s}^2 & = 2 \langle \Lambda^s \partial_t J_{\varepsilon}u, \Lambda^s J_{\varepsilon}u\rangle_{L^2} \nonumber \\
& & = -2 \langle \Lambda^s J_{\varepsilon}(u^2 u_x), \Lambda^s J_{\varepsilon}u\rangle_{L^2} -3 \langle \Lambda^{s-1} J_{\varepsilon}(u^3), \Lambda^s J_{\varepsilon}u\rangle_{L^2} \label{eq:pers1}\\
& & \ \ -3 \langle \Lambda^{s-1} J_{\varepsilon}(u u_x^2), \Lambda^s J_{\varepsilon}u\rangle_{L^2} - \langle \Lambda^{s-2} J_{\varepsilon}(u_x^3), \Lambda^s J_{\varepsilon}u\rangle_{L^2} \label{eq:pers2}.
\end{eqnarray}
We write the first summand in \eqref{eq:pers1} as a sum of two terms 
\begin{equation} \label{eq:pers1ab}
\langle \Lambda^s J_{\varepsilon}(u^2 u_x), \Lambda^s J_{\varepsilon}u\rangle_{L^2}=\langle [\Lambda^s,u^2 \partial_x]u, \Lambda^s J_{\varepsilon}^2 u\rangle_{L^2} 
+ \langle u^2 \partial_x \Lambda^s u, \Lambda^s J_{\varepsilon}^2 u\rangle_{L^2}.
\end{equation}
and estimate each term separately. For the first term we use a commutator estimate that can be found in \cite{Taylor}:
\[ \|[\Lambda^s, u^2\partial_x]u \|_{L^2}\leq \|u^2\|_{H^s}\|u_x\|_{\infty}+\|u^2\|_{C^1}\|u_x\|_{H^{s-1}}
\]
\[ \leq \|u\|_{C^1}^2 \|u\|_{H^s}.
\]
For the second summand on the right hand side of \eqref{eq:pers1ab} we use the identity 
$$\int_{S^1} u^2(\Lambda^s\partial_x u)(\Lambda^s J_{\varepsilon}^2u)dx = -\int_{S^1}uu_x (\Lambda^s J_{\varepsilon}u)^2 dx$$ 
to obtain 
\[ \langle u^2 \partial_x \Lambda^s u, \Lambda^s J_{\varepsilon}^2 u\rangle_{L^2}\lesssim \|u\|_{C^1}^2 \|J_{\varepsilon}u\|_{H^s}^2.
\] 

For the remaining three terms in \eqref{eq:pers1}-\eqref{eq:pers2} we do not need to use mollifiers. They are bounded by $\big(\| \Lambda^{s-1}(u^3+uu_x^2)\|_{L^2}+\| \Lambda^{s-2}(u_x^3)\|_{L^2}\big)\|u\|_{H^s}$ by Cauchy-Schwartz. Therefore, passing to the limit $\varepsilon\rightarrow 0$ yieds the inequatity \eqref{eq:persgoal}.
\proofend

Novikov's equation \eqref{eq:nov2} is a generalized right Euler-Poincar\'e equation
\begin{equation}\label{teta}
\frac{d}{dt}\frac{\delta l}{\delta u}=-\theta^*_u\frac{\delta l}{\delta u},
\end{equation}
where the right invariant Lagrangian is determined by $l=\frac{1}{2}\int mu dx$ and $\theta^*_u =u\partial_x+\frac{3}{2}u'$ determines an infinitesimal action of the Lie algebra of vector fields $\mbox{Vect}(S^1)$ on its dual $\mathfrak g^*$.  This generalization of Euler-Lagrange equations in \cite{TV} is motivated by systems whose configuration spaces are given by Lie group representations other than the coadjoint representation. Indeed the Lie algebra action $\theta^*$ replaces the Lie algebra coadjoint action $\mbox{ad}^*$. The corresponding group action is given by $\Theta^*_\eta (m)=(m\circ\eta)(\partial_x\eta)^{3/2}$ replacing the group coadjoint action $\mbox{Ad}^*$. Then a conserved quantity along solutions $u$ of Novikov's equation is
\begin{equation}\label{eq:orbit}
\Theta^*_\eta(m)=(m\circ\eta)(\partial_x\eta')^{3/2}=u_0-u_0''.
\end{equation}
for any curve $\eta$ in $\mbox{Diff}(S^1)$ satisfying $u=\dot{\eta}\circ\eta^{-1}$ \cite{TV}. We use this conserved quantity in the proof of the global existence and uniqueness theorem that follows.

\begin{theorem}[Global existence and uniqueness]\label{th:gwp}
For $s\mathfrak \geq 3$ assume that $u_0\in H^s({\Bbb T})$ and
\begin{equation}\label{eq:sign_cond}
\Lambda^2 u_0\mathfrak \geq 0 \quad (\mbox{or } \leq 0).
\end{equation}
Then the Cauchy problem \eqref{eq:nloc}-\eqref{eq:data} has a unique global (in time) solution $u\in C({\Bbb R},H^s({\Bbb T}))\cap C^1({\Bbb R}, H^{s-1}({\Bbb T}))$.
\end{theorem}
\proofbegin
First we note that, by proposition \ref{prop:persistence}, it is sufficient to find a time independent bound for $\|u\|_{C^1}$ in order to extend the local (in time) solutions of theorem \ref{th:lwp}. Note also that we have $\| u\|_{C^1}\lesssim \|\Lambda^2 u\|_{L^1}$ by the Sobolev imbedding theorem. Furthermore the orbit invariant in \eqref{eq:orbit} guarantees that $\Lambda^2 u_0\mathfrak \geq 0$ implies $\Lambda^2 u \mathfrak \geq 0$ for all $t$. Therefore we have
\[
\|u\|_{C^1}\lesssim \|\Lambda^2u\|_{L^1}=\int_{S^1}\Lambda^2 u \ dx = \int_{S^1}u \ dx \leq 1+\|u_0\|_{H^1}\]
and hence the unique solution of theorem \ref{th:lwp} exists for all time.
\proofend
\section{Analyticity of solutions}

In this section we look for real analytic solutions of the periodic Cauchy problem \eqref{eq:nloc}-\eqref{eq:data}. The classical Cauchy-Kowalevski theorem does not apply to equation \eqref{eq:nov1}. However a contraction argument on a scale of Banach spaces can be used for the nonlocal form of this equation to prove the following theorem. 

\begin{theorem}\label{th:analytic}
If the initial data $u_0$ is analytic on ${\Bbb T}$ then there exists an $\varepsilon>0$ and a unique solution $u(t,x)$ to the Cauchy problem \eqref{eq:nloc}-\eqref{eq:data} that is analytic both in $x$ and $t$ on ${\Bbb T}$ for all $t$ in $(-\varepsilon,\varepsilon)$.
\end{theorem} 

\begin{remark}
The contraction argument used in the proof of theorem \ref{th:analytic} is on a decreasing scale of Banach spaces.  We give the statement of this abstract theorem here as stated in \cite{Nis} and \cite{BG} for the convenience of the reader:

Given a Cauchy problem 
\begin{equation}
\partial_t u = F(t, u(t)), \quad u(0)=0,
\label{eq:analytic}
\end{equation}
 and a decreasing scale of Banach spaces $\{ X_{s}\}_{0<s<1}$ so that for any $s'<s$ we have $X_{s}\subset
X_{s'}$ and $ ||| . |||_{s'} \leq ||| . |||_{s} $, 
let $T,R$ and $C$ be positive numbers, suppose that $F$ satisfies the following conditions:

1.) If for $0<s'<s<1$ the function $t\longmapsto u(t)$ is holomorphic in $|t|<T$ and continuous on $|t|\leq
T$ with values in $X_{s}$ and $\sup_{|t|\leq T} ||| u(t) |||_{s} <R,$
then $t\longmapsto F(t,u(t))$ is a holomorphic function on $|t|<T$ with values in $X_{s'}$.

2.) For any $0<s'<s\leq 1$ and any $u,v \in X_{s}$ with $|||u|||_{s}<R, ||| v |||_{s}<R$,
\[ \sup_{|t|\leq T} ||| F(t,u)-F(t,v) |||_{s'} \leq \frac{C}{s-s'}|||u-v |||_{s}.
\]

3.) There exists $M>0$ such that for, any $0<s<1$, 
\[ \sup_{|t|\leq T} |||F(t,0) |||_{s} \leq \frac{M}{1-s} .
\]
Then there exists a $T_{0}\in (0,T)$ and a unique function $u(t)$, which for every $s \in (0,1)$ is
holomorphic in $|t|<(1-s)T_{0}$ with values in $X_{s}$, and is a solution to the initial value problem 
(\ref{eq:analytic}).
\end{remark}

Next we restate the Cauchy problem \eqref{eq:nloc}-\eqref{eq:data} in a more convenient form. Let $v=u_x$. Then the problem \eqref{eq:nloc}-\eqref{eq:data} can be written as a system for $u$ and $v$:
\begin{equation}\label{eq:system}
\left\{
\begin{array}{l}
u_t=-u^2v-\Lambda^{-2}(3u^2 v+2v^3+3uvv_x)=F(u,v), \\
v_t=-uv^2-u^2 v_x- \Lambda^{-2}\partial_x (3u^2 v+2v^3+3uvv_x)=G(u,v), \\
u(0,x)=u_0(x), \quad v(0,x)=u_0'(x).
\end{array}
 \right.
 \end{equation}
This is our Cauchy problem as in \eqref{eq:analytic}.

For $s>0$ let $E_s$ be defined as
\[ E_s=\left\{ u\in C^{\infty}({\Bbb T}): ||| u |||_{s}=\sup_{k>0}\frac{\| \partial_x^k u \|_{H^2}s^k}{k!/(k+1)^2}<\infty \right\}.
\]
This norm is introduced in \cite{HM3} to study the analyticity of the Cauchy problem for Camassa-Holm equation. Let $X_s$ denote the product space $E_s\times E_s$ equipped with a product norm $||| \ . \ |||_{X_s}$. We use the decreasing scale of Banach spaces $X_s$ with $s>0$ for the contraction argument.
  
We would like to include three properties of the spaces $E_s$ that we use in the proof of theorem \ref{th:analytic}. The first one can be seen as a ring property for $E_s$ when $0<s<1$: For any $u, v \in E_s$ there is a constant $c>0$ that is independent of $s$ such that
\begin{equation}\label{eq:sproduct}
|||uv|||_s\leq c|||u|||_s|||v|||_s.
\end{equation}  
The second property is an estimate for the differential operator $\partial_x$: For $0<s'<s<1$ we have
\begin{equation}\label{eq:prop_dx}
|||u_x|||_{s'}\leq \frac{C}{s-s'}|||u|||_s
\end{equation} 
for all $u$ in $E_s$. 
The third property is a pair of estimates for the operator $\Lambda^{-2}$: For $0<s<1$ we have 
\begin{equation}\label{eq:prop_lambda}
|||\Lambda^{-2}u|||_{s'}\leq |||u|||_s \quad \mbox{ and } \quad
|||\Lambda^{-2}\partial_x u|||_{s'}\leq |||u|||_s 
\end{equation}
for all $u$ in $E_s$. We refer to \cite{HM3} for detailed proofs of these properties of $E_s$. Now we have all the tools we need to prove theorem \ref{th:analytic}.

\pott 
Note that it is sufficient to verify the conditions 1.) and 2.) in the statement of the abstract Cauchy-Kowalevski theorem above for both $F(u,v)$ and $G(u,v)$ in the system \eqref{eq:system} since neither $F$ nor $G$ depend on $t$ explicitly.    

We observe that, for $0<s'<s<1$, the estimates \eqref{eq:sproduct}, \eqref{eq:prop_dx} and \eqref{eq:prop_lambda} imply the bounds
\[ |||F(u,v)|||_{s'}\leq 4|||u|||_s^2 |||v|||_s + 2|||v|||_s^3 + \frac{C_1}{s-s'}|||v|||_s^2 |||u|||_s
\]
and 
\[ |||G(u,v)|||_{s'}\leq |||u|||_s |||v|||_s^2 + \frac{C_2}{s-s'}|||u|||_s^2 |||v|||_s+ |||v|||_s^3 + \frac{C_3}{s-s'}|||v|||_s^2 |||u|||_s,
\]
hence condition 1.) holds.

Note that to verify the second condition it suffices to estimate $|||F(u_1,v)-F(u_2,v)|||_{s'}$ and $|||F(u,v_1)-F(u,v_2)|||_{s'}$  separately. By the estimates \eqref{eq:sproduct}, \eqref{eq:prop_dx} and \eqref{eq:prop_lambda} we have 
\[ |||F(u_1,v)-F(u_2,v)|||_{s'}\leq C_R \left(|||u_1-u_2|||_s + \frac{1}{s-s'}|||u_1-u_2|||_s \right)
\]
where the constant $C_R$ depends only on $R$. Similarly we have 
\[ |||F(u,v_1)-F(u,v_2)|||_{s'}\leq \bar{C}_R\left( |||v_1-v_2|||_s + \frac{1}{s-s'}|||v_1-v_2|||_s\right). 
\]
Furthermore $ |||G(u_1,v)-G(u_2,v)|||_{s'}$ and $ |||G(u,v_1)-G(u,v_2)|||_{s'}$ are bounded respectively by $C_R|||u_1-u_2|||_s + \bar{C}_R/(s-s')|||u_1-u_2|||_s $ and $C_R|||v_1-v_2|||_s + \bar{C}/(s-s')|||v_1-v_2|||_s$ using the same argument for a new set of constants $C_R$ and $\bar{C}_R$. Therefore condition 2.) holds as well.

\proofend

\bibliography{is}

\end{document}